\documentclass[12pt]{article}
\usepackage{amssymb,amsmath,cite}
\usepackage[dvips]{graphicx,color}
\usepackage{psfrag,subfigure}
\title{Divergent Integrals of Certain Analytical Functions In the Sense of Zeta Regularization}
\author{Farhad Aghili\footnote{email: faghili@encs.concordia.ca $\;$ or $\;$ farhad.aghili@gmail.com}}
\date{}

\begin{document}

\maketitle

\begin{abstract}
In this paper we extend the Zeta function regularization technique, which gives a meaningful solution to  divergent power series, in order to assign finite values to divergent integral of certain transcendental   functions $f(x)$. The functions are assumed to be analytical and hence they have a convergent Maclaurin series with infinity radius of convergence. Using Maclaurin series and binomial expansion, we  equivalently  convert the divergent integral to an infinite series in terms of Riemann zeta function. It is shown that
the  infinite series is convergent under a certain condition and consequently the solution of the divergent integral becomes $\int_0^{\infty} f(x) dx = \sum_{k=0}^{\infty} \frac{(-1)^{k+1}}{(k+2)!} f^{(k)}(0)$.  The advantage of this Zeta regularization technique is that it does not not require introducing any cut-off function to calculate such divergent integrals. This formula is applied to calculate the divergent integral of several common transcendental functions. 
\end{abstract}

\section{Introduction}
Divergent series and integrals appear in several areas of mathematics and modern physics. 
Zeta function regularization  method assigns finite values to divergent sums  and this technique is commonly applied  in number theory and  theoretical physics to give precise meanings to ill-conditioned sums. Zeta function regularization is conceptualized from one of the most  enigmatic results in arithmetic that the sum of an infinite divergent series can result in a finite number which is mathematically consistent. In spite of being a peculiar and counter-intuitive abstraction in mathematics, divergent series appear also in physics where experimental observation matches the prediction of the theoretical results~\cite{Dowling-1989,Schumayer-Hutcheiston-2011}.

Solving divergent integral of certain functions is quintessential in mathematics of the quantum field theory (QFT). Introducing suitable cut-off function is the cornerstone of the mathematics of divergent integrals appeared in the calculations of QFT. The definition of the Zeta regularization of a series is extended in \cite{Moreta2-2013}  for calculation of  certain divergent integrals  using the Zeta regularization technique combined with the Euler-Maclaurin summation formula. The Euler-Maclaurin formula is extended into a recursive formula to provide finite regularization to divergent integrals.  
An explicit bounded solution of divergent integral of power functions was presented in \cite{Aghili-Tafazoli-2018}. This work extends the results of \cite{Aghili-Tafazoli-2018} for the case of divergent integral of certain analytical functions. Using the notion of Zeta function regularization, we will  show that the divergent integrals can be represented by a infinite series, which is convergent under certain condition.  Subsequently, the formula is  applied to some common divergent integrals to derive  finite solutions in the sense of Zeta regularization.

\section{Zeta Regularization of Divergent Integrals}
Suppose the following integral of function $f(x)$ diverges  
\begin{equation} \label{eq:div_integral}
I=\int_{0}^{\infty} f(x) dx
\end{equation}
Also assume that the function has a convergent Maclaurin series with infinity radius of convergence, i.e.,  
\begin{equation} \label{eq:Tylor}
f(x) = \sum_{k=0}^{\infty} \frac{x^k}{k!} f^{(k)}(0)    \qquad \forall x \in \mathbb{R}
\end{equation}
Then we seek to assign a bounded meaningful solution to the divergent integral \eqref{eq:div_integral} using a Zeta regularization technique. Note that if function $f(z)$ is equal to its Maclaurin series for all $z$ in the complex plane, then it is called entire. Examples of entire functions include polynomials, exponential function, and the trigonometric functions sine and cosine, while square root and the fractional functions are not entire. It is worthwhile mentioning that the necessary and sufficient for function $f(z)$ in terms of its complex variable $z=x+iy$ to be analytic are that: i) $\partial u/\partial x$, $\partial u/\partial y$,  $\partial v/\partial x$, $\partial v/\partial y$ are continuous and ii) the aforementioned  four partial derivatives of the real part $u(x, y)$ and imaginary 
part $v(x, y)$ of the function satisfies the Cauchy-Riemann equations
\begin{equation}
\frac{\partial u}{\partial x}= \frac{\partial v}{\partial y} \quad \mbox{and} \quad \frac{\partial u}{\partial y}= -\frac{\partial v}{\partial x}
\end{equation}

Replacing the function with its Maclaurin series \eqref{eq:Tylor} in the integral \eqref{eq:div_integral} and then switching the integral and sum, we can rewrite the divergent integral as follows
\begin{equation} \label{eq:split}
I= \sum_{k=0}^{\infty} \frac{f^{(k)}(0)}{k!} \int_0^{\infty}x^k \; dx 
\end{equation}
The above improper integral can be equivalently written as the following infinite series by splitting the range of integration limits into successive integer numbers
\begin{equation} \label{eq:two_integrals}
\int_0^{\infty}x^k \; dx = \sum_{n=1}^{\infty} \int_{n-1}^n x^k dx = \sum_{n=1}^{\infty} \frac{1}{k+1}\Big( n^{k+1} - (n-1)^{k+1}  \Big)
\end{equation}  
Using \eqref{eq:two_integrals} in \eqref{eq:split}  yields
\begin{equation} \label{eq:I_n+1}
I= \sum_{r=0}^{\infty} \frac{f^{(k)}(0)}{(k+1)!}  \sum_{n=1}^{\infty} \Big( n^{k+1} - (n-1)^{k+1}  \Big)
\end{equation}
According to the binomial  polynomial expansion of $(n-1)^{k+1}$, we have  
\begin{subequations}
\begin{equation} \label{eq:n-1_binomial}
(n-1)^{k+1} = \sum_{p=0}^{k+1}  {{k+1}\choose{p}} (-1)^{k+1-p}n^p,
\end{equation}
where the binomial coefficient is the defined as
\begin{equation} \label{eq:choose}
{{m}\choose{l}} = \frac{m !}{l !(m-l) !} 
\end{equation}
\end{subequations}
Upon substitution of \eqref{eq:n-1_binomial} in \eqref{eq:I_n+1}, we arrive at
\begin{align} \label{eq:}
I &= \sum_{r=0}^{\infty} \frac{f^{(k)}(0)}{(k+1)!} \sum_{n=1}^{\infty} \sum_{p=0}^{k} {{k+1}\choose{p}} (-1)^{r-p}n^p \\ \label{eq:Delta_r2} 
& = \sum_{k=0}^{\infty} \frac{f^{(k)}(0)}{(k+1)!} \sum_{p=0}^{k} {{k+1}\choose{p}} (-1)^{k-p}\sum_{n=1}^{\infty} n^p \\
& = \sum_{k=0}^{\infty} \frac{f^{(k)}(0)}{(k+1)!} \sum_{p=0}^{k} {{k+1}\choose{p}} (-1)^{k-p}\zeta(-p) 
\end{align}
For positive integer $p\geq 0$, the Zeta function is related to Bernoulli numbers by
\begin{equation} \label{eq:zeta_B}
\zeta(-p) = (-1)^p \frac{B_{p+1}}{p+1}
\end{equation}
Thus, we have
\begin{equation}
I= \sum_{k=0}^{\infty} \frac{(-1)^{k} f^{(k)}(0)}{(k+1)!} \sum_{p=0}^{k} {{k+1}\choose{p}} \frac{B_{p+1}}{p+1}
\end{equation}
Moreover, from definition \eqref{eq:choose}, one can verify that the successive binomial coefficients hold the following useful identity   
\begin{equation} \label{eq:2Cs}
{{k+1}\choose{p}}  = \frac{p+1}{k+2} {{k+2}\choose{p+1}}
\end{equation}
Thus we have
\begin{align} \notag
I &= \sum_{k=0}^{\infty} \frac{(-1)^{k} f^{(k)}(0)}{(k+2)!} \sum_{p=0}^{k} {{k+2}\choose{p+1}} B_{p+1}\\ \notag
&= \sum_{k=0}^{\infty} \frac{(-1)^{k} f^{(k)}(0)}{(k+2)!} \sum_{p=1}^{k+1} {{k+2}\choose{p}} B_{p}\\ \notag
&= \sum_{k=0}^{\infty} \frac{(-1)^{k} f^{(k)}(0)}{(k+2)!} \left( \sum_{p=0}^{k+1} {{k+2}\choose{p}} B_{p} - {{k+2}\choose{0}} B_{0} \right) \\ \label{eq:I_sum2}
&= \sum_{k=0}^{\infty} \frac{(-1)^{k} f^{(k)}(0)}{(k+2)!} \left( \sum_{p=0}^{k+1} {{k+2}\choose{p}} B_{p} - 1 \right) 
\end{align}
On the other hand, the Bernoulli numbers satisfy the following property 
\begin{equation} \label{eq:sumCB}
\sum_{p=0}^{m-1} {{m}\choose{p}}  B_p = \sum_{p=0}^{k+1} {{k+2}\choose{p}} B_p = 0
\end{equation}
Finally using identity  \eqref{eq:sumCB} in \eqref{eq:I_sum2}, we arrive at the explicit expression of the divergent integral based on Zeta regularization   
\begin{equation} \label{eq:divergent_integral}
\boxed{\int_0 ^{\infty} f(x) \; dx = \sum_{k=0}^{\infty} \frac{(-1)^{k+1} }{(k+2)!}f^{(k)}(0)}
\end{equation}
Invoking the ratio test, one can conclude that that the above infinite series converge if
\begin{equation}
\lim_{k \rightarrow \infty}   \frac{\frac{(-1)^{k+2} f^{(k+1)}(0)}{(k+3)!}} {\frac{(-1)^{k+1} f^{(k)}(0)}{(k+2)!}}<1
\end{equation}
In other words, the convergence of the above series is
\begin{equation}
\lim_{k \rightarrow \infty}  \frac{|f^{(k+1)}(0)|}{(k+3)|f^{(k)}(0)| } <1
\end{equation}

Equation \eqref{eq:divergent_integral} can be used to derive Zeta regularization of many divergent integrals. For instance, if the  analytical function is selected to be a polynomial 
\begin{equation}
f(x) = x^m
\end{equation}
Then,
\begin{equation} \label{eq:fr_polynomial}
f^{(k)}(0) = \left\{ \begin{array}{ll} m! \quad & \mbox{if} \;\; k=m\\ 0 & \mbox{otherwise} \end{array} \right.
\end{equation} 
Applying \eqref{eq:divergent_integral} to \eqref{eq:fr_polynomial} yields
\begin{equation} \notag
\int_0^{\infty} x^m \; dx  = 0+0+\cdots + \frac{(-1)^{m+1} }{(m+2)!} m! + 0 + 0+ \cdots,
\end{equation}
and thus
\begin{equation} \notag
\int_0^{\infty} x^m \; dx =  \frac{(-1)^{m+1} }{(m+1)(m+2)},
\end{equation}
which is consistent with the previous result in \cite{Aghili-Tafazoli-2018}.
 
\section{Examples of Divergent integrals}
Consider the trigonometric function
\begin{equation}
f(x) =\sin(x)
\end{equation}
whose derivatives are $f^{(2k+1)}(0)=(-1)^{2k+1}$. Thus, according to \eqref{eq:divergent_integral} we have
\begin{equation}
\int_0^{\infty} \sin(x) \; dx = \sum_{k=0}^{\infty} \frac{(-1)^{k+1}}{(2k+3)!}
\end{equation}
By inspection one can show that the right-hand side of the above equation is equal to $\sin(1)-1$. 

Examples of  Zeta regularization of divergent integrals of several  transcendental functions  are given below:
\begin{subequations}
\begin{align}
&\int_0^{\infty} \sin(x) dx = \sum_{k=0}^{\infty} \frac{(-1)^{k+1}}{(2k +3)!} =\sin(1)-1 \\
&\int_0^{\infty} \cos(x) dx = \sum_{k=0}^{\infty} \frac{(-1)^{2k+1}}{(2k +2)!} =\cos(1)-1 \\
&\int_0^{\infty} \sinh(x) dx = \sum_{k=0}^{\infty} \frac{1}{(2k +3)!} =\sinh(1)-1 \\
&\int_0^{\infty} \cosh(x) dx = \sum_{k=0}^{\infty} \frac{(-1)^{2k+1}}{(2k +2)!} =1-\cosh(1) \\
&\int_0^{\infty} e^x dx = \sum_{k=0}^{\infty} \frac{(-1)^{k+1}}{(k +2)!} = -\frac{1}{e} \\
& \int_0^{\infty} \log(1+x) dx = \sum_{k=0}^{\infty} \frac{1}{k(k+1)(k+2)} = \frac{1}{4}\\
&\int_0^{\infty} e^{x^2} dx = \sum_{k=0}^{\infty} \frac{(-1)^{k+1}}{2(k +2)(k+1)k!} = \frac{1}{2} - \frac{1}{2e} - \frac{\sqrt{\pi}}{e} \mbox{erf}(1) \\
&\int_0^{\infty} \mbox{erf}(x)dx = \frac{2}{\sqrt{\pi}} \sum_{k=0}^{\infty} \frac{(-1)^k}{2(2k+1)(k+1)(2k+3)k!} = \frac{3}{2} \mbox{erf}(1) + \frac{1}{2e \sqrt{\pi}} - \frac{1}{\sqrt{\pi}}\\
&\int_0^{\infty} e^{-x^2}dx = \frac{2}{\sqrt{\pi}} \sum_{k=0}^{\infty}
\frac{(-1)^{k+1}}{2(2k+1)(k+1)!} = -\frac{\sqrt{\pi}}{2} \mbox{erf}(1)+\frac{1-e}{2e}\\
&\int_0^{\infty} x \sin(x)dx =  \sum_{k=0}^{\infty} \frac{1}{2(k+2)(2k+3)(1+2k)!} = 2 + \sinh(1) -2 \cosh(1)
\end{align}
\end{subequations}

\section{Conclusions}

We extended the zeta function regularization technique in order to assign finite values to divergent integral of a class of analytical functions which satisfy a convergence condition. Using the binomial expansion and Maclaurin series, it has been shown that such divergent integrals can be expressed by infinite series in terms of the Riemann zeta function. Subsequently, it has been shown that the 
divergent integral can be written into the following infinite series $\int_0^{\infty} f(t) dt = \sum_{k=0}^{\infty} \frac{(-1)^{k+1}}{(k+2)!} f^{(k)}(0)$, which becomes convergent under certain condition. Calculation of several common divergent integrals using this formula has been presented.

\bibliographystyle{IEEEtran}

\end{document}